\documentstyle[leqno,draft]{article}
\begin{document}
\renewcommand{\thesection}{\arabic{section}.}
\renewcommand{\theequation}{\arabic{section}.\arabic{equation}}
\newcommand{\be}{\begin{eqnarray}}
\newcommand{\en}{\end{eqnarray}}
\newcommand{\no}{\nonumber}
\newcommand{\la}{\lambda}
\newcommand{\laa}{\Lambda}
\newcommand{\ep}{\epsilon}
\newcommand{\de}{\delta}
\newcommand{\D}{\Delta}
\newcommand{\pl}{\parallel}
\newcommand{\ov}{\overline}
\newcommand{\bet}{\beta}
\newcommand{\al}{\alpha}
\newcommand{\fr}{\frac}
\newcommand{\pa}{\partial}
\newcommand{\we}{\wedge}
\newcommand{\om}{\Omega}
\newcommand{\na}{\nabla}
\newcommand{\lan}{\langle}
\newcommand{\ra}{\rangle}
\newcommand{\fa}{\sum_{\alpha=1}^n}
\newcommand{\vs}{\vskip0.3cm}
\renewcommand{\thefootnote}{}
\title {Inequalities for eigenvalues of \\ the buckling problem of arbitrary  order }

\author{Qing-Ming Cheng, \  Xuerong Qi, \  Qiaoling Wang and  Changyu Xia}

\date{} 
\maketitle
\begin{abstract} \noindent
This paper studies  eigenvalues of the buckling problem of arbitrary
order on bounded domains  in  Euclidean spaces and spheres. We prove
universal bounds for the $k$-th eigenvalue in terms of the lower  ones independent  of  the domains.
Our results strengthen the recent work in \cite{[JLWX1]} and generalize Cheng-Yang's recent estimates \cite{[CY8]} on the buckling eigenvalues of order two to arbitrary order.
\end{abstract}


\section{Introduction}
Let $\om$ be a bounded domain with smooth boundary in an  $n(\geq 2)$-dimensional Riemannian manifold $M$  and denote by $\Delta $ the Laplace operator acting on functions on $M$. Let $\nu$ be the outward unit normal vector field of $\pa \om$ and
 let us consider the following  eigenvalue problems :
\be& &
~\Delta u=-\la u \ \ \ \ \ {\rm in \ \ } \om,  \ \ \ \
u=0, \ \ \ \ \ \ \ \ \ \ {\rm on \ \  } \pa \om,
\\  & &
\Delta^2 u= -\laa \Delta u \ \ {\rm in \ \ } \om,  \ \ \ \
u=\fr{\pa u}{\pa \nu }=0, \ \  {\rm on \ \ } \pa \om.
\en
They are called the {\it  fixed membrane problem} and
the {\it bucking problem}, respectively. It should be mentioned that the buckling problem (1.2) has interpretations in physics, that is, it  describes the critical buckling load of a clamped plate subjected to a uniform
compressive force around its boundary.
Let
\be\no
& & 0<\la_1<\la_2\leq\la_3\leq\cdots,\\
\no & & 0<\laa_1\leq\laa_2\leq\laa_3\leq\cdots
\en
denote the successive eigenvalues for (1.1) and (1.2), respectively.
Here each eigenvalue is repeated according to its
multiplicity. An
important theme of geometric analysis is to estimate these (and other) eigenvalues. When $\om$ is a bounded domain in an
$n$-dimensional Euclidean space ${\bf R}^n$,
Payne, P\'olya and Weinberger (cf. \cite{[PPW1]},\cite{[PPW2]}) proved the bound
\be
\la_{k+1}-\la_k\leq \fr 4{kn}\sum_{i=1}^k\la_i, \ \ k=1, 2,\cdots.
\en
Inequality of this type is called a universal inequality since it does not depend on $\om$.

On the other hand,  Payne, P\'olya and Weinberger also studied
eigenvalues of  the buckling problem (1.2) for  bounded domains in ${\bf R}^n$  and
proved (cf. \cite{[PPW1]},\cite{[PPW2]})
$$\laa_2/\laa_1<3 \ \ \ \ {\rm for }\ \om\subset {\bf R}^2.
$$
For $\om\subset {\bf R}^n$,  this reads
$$
\laa_2/\laa_1<1+4/n.
$$
Furthermore, Payne, P\'olya and Weinberger  proposed  the following
\vskip0.3cm
\noindent
{\bf Problem 1 } (cf. \cite{[PPW1]},\cite{[PPW2]}). {\it  Can one obtain a universal inequality for the
  eigenvalues of the buckling problem {\rm (1.2)} on a bounded domain in ${\bf R}^n$
  which is similar to the universal inequality {\rm (1.3)} for the eigenvalues of the fixed membrane problem {\rm (1.1)}  ?
}
\vskip0.3cm
 With respect to the above problem,
 Hile and Yeh \cite{[HY]} obtained
$$
\fr{\laa_2}{\laa_1}\leq \fr{n^2+8n+20}{(n+2)^2}\ \ \ \ \ {\rm for \ } \om\subset {\bf R}^n.
$$
Ashbaugh \cite{[As1]} proved :
\be
\sum_{i=1}^n\laa_{i+1}\leq(n+4)\laa_1.
\en
This inequality has been improved  to the following form in \cite{[JLWX2]}:
\be \sum_{i=1}^n\laa_{i+1}
+\fr{4(\laa_2-\laa_1)}{n+4}\leq (n+4)\laa_1.
\en
By introducing a new method of constructing trial functions, Cheng and Yang  \cite{[CY3]} obtained the following
universal inequality  and thus solved the above problem:
\be \sum_{i=1}^k(\laa_{k+1}-\laa_i)^2\leq
\fr{4(n+2)}{n^2}\sum_{i=1}^k(\laa_{k+1}-\laa_i)\laa_i.
\en
Recently, Cheng-Yang \cite{[CY8]} have  proved the following inequality:
\be n\sum_{i=1}^k(\laa_{k+1}-\laa_i)^2\leq \left(n+\fr 43\right)\sum_{i=1}^k \delta_i (\laa_{k+1}-\laa_i)^2+\sum_{i=1}^k \fr 1{\delta_i} (\laa_{k+1}-\laa_i)^2\laa_i,
\en
where $\{\delta_i\}_{i=1}^k$ is any positive non-increasing monotone sequence.
Taking
\be\no
\delta_i=\sqrt{\fr{\sum_{i=1}^k (\laa_{k+1}-\laa_i)^2\laa_i}{\left(n+\fr 43\right)\sum_{i=1}^k (\laa_{k+1}-\laa_i)^2}},\ i=1,\cdots, k
\en
in (1.7), Cheng-Yang have obtained
\be \sum_{i=1}^k(\laa_{k+1}-\laa_i)^2\leq \fr{4\left(n+\fr 43\right)}{n^2}\sum_{i=1}^k(\laa_{k+1}-\laa_i)\laa_i,
\en
which is stronger than (1.6).

It has been proved in \cite{[WX2]} that for the problem (1.2) if $\om$ is a  domain in an $n$-dimensional unit sphere  $S^n$, then we have
\be & &
2\sum_{i=1}^k(\laa_{k+1}-\laa_i)^2\\ \no &\leq&
\sum_{i=1}^k(\laa_{k+1}-\laa_i)^2\left(\delta \laa_i+\fr{\delta^2(\laa_i-(n-2))}{4(\delta\laa_i+n-2)}\right)  +\fr 1{\delta}\sum_{i=1}^k(\laa_{k+1}-\laa_i)\left(\laa_i+\fr{(n-2)^2}4\right),
\en
where $\delta $ is any positive constant. This inequality has been improved recently in \cite{[CY8]}  and \cite{[JLWX1]}, respectively.
\vskip0.3cm
In this paper, we will investigate  eigenvalues of the  buckling problem of arbitrary  order:
\be\left\{\begin{array}{cl}
 (-\Delta)^l u= -\laa\D u, \ \ \ \ \ \ \ \ \ \ {\rm in} \ \ \om, \\
   u=\fr{\pa u}{\pa \nu}=\cdots =\fr{\pa^{l-1} u}{\pa \nu^{l-1}}=0,\ {\rm on} \ \ \pa \om,
\end{array}\right.\en
where $ \om$ is a  bounded domain in a Euclidean space or a unit sphere and $l$ is any integer no less than $2$.
Yang type inequalities for  eigenvalues of the problem (1.10) have been obtained recently in \cite{[JLWX1]}.  In this paper, we  prove :
\vskip0.3cm
\noindent
{\bf Theorem 1.1.}
{\it Let $\laa_i$ be the $i$-th eigenvalue of the buckling problem {\rm (1.10)}, where $\om$ is a  bounded domain with smooth boundary in  ${\bf R}^n$. Then for any positive non-increasing monotone sequence
$\{\delta_i\}_{i=1}^k$,
we have
\be & &
n\sum_{i=1}^k(\laa_{k+1}-\laa_i)^2\\ \no &\leq& 
~\sum_{i=1}^k \delta_i(\laa_{k+1}-\laa_i)^2\left(2l^2+\left(n-\fr{14}3\right)l+\fr 83-n\right)\laa_i^{(l-2)/(l-1)}
+ \sum_{i=1}^k \fr 1{\delta_i} (\laa_{k+1}-\laa_i)\laa_i^{1/(l-1)}.
\en
}
\vskip0.3cm
\noindent
{\bf Remark 1.1.} When $l=2$, (1.11) becomes Cheng-Yang's inequality (1.7).
\vskip0.3cm

\noindent
{\bf Remark 1.2.} Taking
\be\no
\de_1=\de_2=\cdots \de_k=\left\{ \fr {\sum_{i=1}^k (\laa_{k+1}-\laa_i)\laa_i^{1/(l-1)}}{ \left(2l^2+\left(n-\fr{14}3\right)l+\fr 83-n\right)\sum_{i=1}^k(\laa_{k+1}-\laa_i)^2\laa_i^{(l-2)/(l-1)}}\right\}^{1/2}
\en
in (1.11), we have
\be
& & \sum_{i=1}^k(\laa_{k+1}-\laa_i)^2\\
\no &\leq& \fr{ 2\left(2l^2+\left(n-\fr{14}3\right)l+\fr 83-n\right)^{1/2} } n\left\{ \sum_{i=1}^k  (\laa_{k+1}-\laa_i)^2\laa_i^{(l-2)/(l-1)} \right\}^{1/2} \left\{\sum_{i=1}^k (\laa_{k+1}-\laa_i)\laa_i^{1/(l-1)}\right\}^{1/2},
\en
which improves the inequality (1.13) in \cite{[JLWX1]}.  From (1.12), we can obtain a   quadratic   inequality  about   $\laa_1,\cdots,\laa_{k+1}$.

\vskip0.3cm 
\noindent
{\bf    Corollary 1.1.} {\it    For    any $k\geq 1$,  the first   $k+1$   eigenvalues of  the buckling problem {\rm (1.10)} with $\om\subset {\bf R}^n$  satisfy the
following   inequality
\be
\sum_{i=1}^k(\laa_{k+1}-\laa_i)^2\leq\fr{ 4\left(2l^2+\left(n-\fr{14}3\right)l+\fr 83-n\right)}{n^2} \sum_{i=1}^k  (\laa_{k+1}-\laa_i)\laa_i.
\en}

\vskip0.3cm
Furthermore, we  prove the following universal inequality for eigenvalues of the buckling problem of arbitrary order on spherical domains.
\vskip0.3cm
\noindent
{\bf Theorem 1.2.} {\it Let $l\geq 2$ and  let $\laa_i$ be the $i$-th eigenvalue of the buckling problem:
\be\left\{\begin{array}{cl}
 (-\Delta)^l u= -\laa\D u, \ \ \ \ \ \ \ \ \ \ {\rm in} \ \ \om, \\
   u= \fr{\pa u}{\pa \nu}=\cdots =\fr{\pa^{l-1} u}{\pa \nu^{l-1}}=0,\ {\rm on} \ \ \pa \om,
\end{array}\right.\en
where $\om$ is a domain  with smooth boundary  in $S^n$. For each $q= 1, \cdots, $ define the polynomials $\Phi_q$  inductively by
\be \Phi_1(t)=t-1, \ \Phi_2(t)=t^2-(n+5)t-(n-2),\en
\be& &
\Phi_q(t)=(2t-2)\Phi_{q-1}(t)-(t^2+2t-n(n-2))\Phi_{q-2}(t),  \ \ q=3,\cdots.
\en
Set
\be
\Phi_{l-1}(t)= t^{l-1}-a_{l-2}t^{l-2}+\cdots +(-1)^{l-2}a_{1}t-(n-2)^{l-2}.
\en
Then for any  positive integer $k$ and any positive non-increasing monotone sequence   $\{\delta_i\}_{i=1}^k$,
we have
\be  & &
\sum_{i=1}^k(\laa_{k+1}-\laa_i)^2\left(2+\fr{n-2}{\laa_i^{1/(l-1)}-(n-2)}\right)
\\ \no &\leq &
 \sum_{i=1}^k (\laa_{k+1}-\laa_i)^2\delta_{i}S_i+
\sum_{i=1}^k\fr{(\laa_{k+1}-\laa_i)}{\delta_i} \left(\laa_i^{1/(l-1)}+\fr{(n-2)^2}4\right),
\en
where
\be
S_i=   \laa_i \left( 1
-\fr 1{\laa_i^{1/(l-1)}-(n-2)}\right)+(-1)^l(n-2)^{l-2}+\sum_{j=1}^{l-2}a_j^{+}\laa_i^{j/(l-1)},
\en
with $a_j^{+}=\max\{a_j, 0\}$ and when $l=2$ we use the convention that $\sum_{j=1}^{l-2}a_j^{+}\laa_i^{j/(l-1)}=0$.
}
\vskip0.3cm
\noindent
{\bf Remark 1.3.} When $l=2$, (1.18) is stronger  than (1.9) and it has been  proved by Cheng-Yang in \cite{[CY8]}.
\vskip0.3cm
\noindent
{\bf Remark 1.4.}
Universal   inequalities   for eigenvalues of  various elliptic    operators  have been    studied  extensively in  recent  years. For the developments in this direction,  we refer to [1-23], [25-30], [33-40] and the references  therein.
\markboth{\hfill  Qiaoling Wang, Changyu Xia \hfill}
{\hfill  Inequalities  for Eigenvalues of the Buckling Problem of Higher Order
 \hfill}
 \vskip0.3cm
\noindent
{\bf Acknowledgement.} The research of the first author is partially supported by a Grant-in-Aid for Scientific Research from JSPS, the research of 
the third author is partially supported by CNPq and the main part of this paper  has been done while the fourth author 
visited to Department of Mathematics, Saga University as a fellowship of JSPS. This author would like to express his 
gratitude to JSPS for finance support and  to Professor Qing-Ming Cheng and Saga University for the worm hospitality.

\section{Proofs of the Results }
\setcounter{equation}{0}
First we recall a method of constructing trial functions developed by Cheng-Yang (cf. \cite{[CY3]}, \cite{[JLWX1]},  
\cite{[WX2]}). Let $M$ be an
  $n$-dimensional complete submanifold in ${\bf R}^{m}$.
Denote by $\langle , \rangle$ the canonical metric on ${\bf R}^{m}$ as well as that induced on $M$. Denote by  $\Delta $ and $\nabla$    the Laplacian and the gradient operator of $M$, respectively. Let $\om$ be a bounded  domain with smooth boundary in $M$  and let $\nu $ be the outward unit normal vector field of $\pa\om$.
For functions $f$ and $g$ on $\om$, the {\it Dirichlet inner product $(f, g)_D$} of $f$ and $g$ is given by
\be \no (f, g)_D=\int_{\om}\langle\nabla f, \ \nabla g\rangle.
\en
The Dirichlet norm of a function $f$ is defined by
\be\no
||f||_D=\{(f, f)_D\}^{1/2}=\left(\int_{\om}|\nabla f|^2\right)^{1/2}.
\en
Consider the
eigenvalue problem
\be\left\{\begin{array}{cl}
 (-\Delta)^l u= -\laa\D u, \ \ \ \ \ \ \ \ \ \ {\rm in} \ \ \om, \\
   u= \fr{\pa u}{\pa \nu}=\cdots =\fr{\pa^{l-1} u}{\pa \nu^{l-1}}=0,\ {\rm on} \ \ \pa \om.
\end{array}\right.\en
Let
\be\no
 0<\laa_1\leq \laa_2\leq\laa_3\leq\cdots
\en denote the successive eigenvalues, where
each eigenvalue is repeated according to its multiplicity.

Let $u_i$ be the $i$-th orthonormal eigenfunction of the problem (2.1) corresponding to the eigenvalue $\laa_i$, $i=1, 2, \cdots,$ that is,
\be\left\{\begin{array}{lll}
 (-\Delta)^l u_i= -\laa_i\D u_i, \ \ {\rm in} \ \ \om, \\
   u_i= \fr{\pa u_i}{\pa \nu}=\cdots =\fr{\pa^{l-1} u_i}{\pa \nu^{l-1}}=0, \ \ {\rm on} \ \ \partial\om,\\
  (u_i, u_j)_D =\int_{\om}\langle \nabla u_i, \nabla u_j\rangle=\delta_{ij}, \ \forall\ i, j.
\end{array}\right.
\en
For $k=1, \cdots, l$, let $\nabla^k$ denote the $k$-th covariant
derivative operator on $M$, defined in the usual weak sense.
 For a function $f$ on $\om$, the squared norm of  $\nabla^k f$ is defined as (cf. \cite{[He]})
\be
\left|\nabla^kf\right|^2=\sum_{i_1,\cdots, i_k=1}^n\left(\nabla^kf(e_{i_1},\cdots, e_{i_k})\right)^2,
\en
where $e_1,\cdots, e_n$ are orthonormal vector fields locally defined on $\om$.
Define the Sobolev space $H_l^2(\om)$ by
$$H_l^2(\om)=\{ f:\ f, \ |\nabla f|,\cdots, \left|\nabla^l f\right|\in L^2(\om)\}.
$$
Then $H_l^2(\om)$ is a Hilbert space with respect to the norm $||\cdot||_{l, 2}$:
\be
||f||_{l, 2}=\left(\int_{\om}\left(\sum_{k=0}^l|\nabla^k f|^2\right)\right)^{1/2}.
\en
Consider the subspace $H_{l,D}^2(\om)$ of $H_l^2(\om)$ defined by
$$H_{l,D}^2(\om)=\left\{f\in H_l^2(\om): \ f|_{\pa \om}=\left. \fr{\pa f}{\pa \nu}\right|_{\pa \om}=\cdots\left. \fr{\pa^{l-1} f}{\pa \nu^{l-1}}\right|_{\pa \om}=0\right\}.
$$
The operator $(-\Delta)^l$ defines a self-adjoint operator acting on $H_{l,D}^2(\om)$
with discrete eigenvalues $0<\laa_1\leq\cdots\leq \laa_k\leq\cdots$ for the buckling problem (2.1) and the eigenfunctions $\{u_i\}_{i=1}^{\infty}$ defined in (2.2)
form a complete orthonormal basis for the Hilbert space $H_{l,D}^2(\om)$. If $\phi\in H_{l,D}^2(\om)$ satisfies $(\phi , u_j)_D=0, \ \forall  j=1, 2, \cdots, k$, then the Rayleigh-Ritz inequality tells us that
\be
\laa_{k+1}||\phi ||_D^2\leq \int_{\om} \phi(-\Delta)^l\phi.
\en
For vector-valued functions $F=(f_1, f_2, \cdots, f_{m}), \ G=(g_1, g_2, \cdots, g_{m}):
\om\rightarrow {\bf R}^{m}$, we define an inner product $(F, G)$ by
$$(F, G)\equiv \int_{\om} \langle F, G\rangle =\int_{\om} \sum_{\alpha =1}^{m} f_{\alpha}g_{\alpha}.$$
The norm of $F$ is given by
$$||F||=(F, F)^{1/2}=\left\{\int_{\om}\sum_{\alpha=1}^{m}f_{\alpha}^2\right\}^{1/2}.$$
Let ${\bf H}_1^2(\om)$ be the Hilbert space of vector-valued functions given by
$${\bf H}_1^2(\om)=\left\{ F=(f_1,\cdots, f_{m}): \om\rightarrow {\bf R}^{m};\  f_{\alpha}, \ |\nabla f_{\alpha}|\in L^2(\om), \ {\rm for} \ \alpha=1,\cdots, m\right\}
$$
with norm
$$||F||_1=\left(||F||^2+\int_{\om}\sum_{\alpha=1}^{m}|\nabla f_{\alpha}|^2\right)^{1/2}.
$$
Observe that a vector field on $\om$ can be regarded as a vector-valued function from $\om$ to ${\bf R}^m$. Let ${\bf H}_{1, D}^2(\om)\subset  {\bf H}_{1}^2(\om)$ be a subspace of ${\bf H}_1^2(\om)$
spanned by the vector-valued functions $\{ \nabla u_i\}_{i=1}^{\infty}$, which form a complete orthonormal basis of ${\bf H}_{1, D}^2(\om)$. For any $f\in H_{l,D}^2(\om), $ we have $\nabla f\in {\bf H}_{1, D}^2(\om)$ and for any $X\in {\bf H}_{1, D}^2(\om)$, there exists a function $f\in H_{l,D}^2(\om)$ such that $X=\nabla f$.

\vskip0.3cm
\noindent
{\bf Lemma 2.1.} (cf. \cite{[JLWX1]},\cite{[JLWX2]}) {\it Let $u_i$ and $\laa_i, i=1, 2, \cdots , $ be as in {\rm (2.2)}, then
\be
0\leq \int_{\om} u_i(-\Delta )^k u_i\leq \laa_i^{(k-1)/(l-1)}, \ \ k=1,\cdots, l-1.
\en
}
\vs
We are now ready to prove the main results in this paper.
\vskip0.3cm
\noindent
{\it Proof of Theorem 1.1.} With the notations as above, we consider now the special case that $\om$ is a  bounded domain in
${\bf R}^{n}$. Denote by $x_1,\cdots, x_n$ the coordinate functions of ${\bf R}^{n}$ and let us
 decompose the vector-valued functions $x_{\alpha}\nabla u_i$ as
\be
x_{\alpha}\nabla u_i=\nabla h_{\alpha i}+ W_{\alpha i},
\en
where
$h_{\alpha i}\in H_{l,D}^2(\om),$ $\nabla h_{\alpha i}$ is the projection of $x_{\alpha} \nabla u_i$ in ${\bf H}_{1, D}^2(\om)$ and $W_{\alpha i}\ \bot\ {\bf H}_{1, D}^2(\om)$. Thus we have
\be
W_{\alpha i}|_{\pa\om}\ =\ 0, \ \ {\rm and} \ \ (W_{\alpha i}, \nabla u)=\int_{\om}
\langle W_{\alpha i}, \nabla u\rangle=0, \ \ {\rm for\ any} \ \ u\in H_{l,D}^2(\om)
\en
and from the discussions in \cite{[CY3]} and \cite{[WX2]} we know that
\be
{\rm div}\ W_{\alpha i}=0,
\en
where for a vector field $Z$ on $\om$, ${\rm div}\ Z$ denotes the divergence of $Z$.

For each $\alpha=1,\cdots, n$, $i=1,\cdots, k$, consider the  functions $\phi_{\alpha i}: \om\rightarrow {\bf R}$, given by
\be
 \phi_{\alpha i}=h_{\alpha i}-\sum_{j=1}^ka_{\alpha ij}u_j,
 \en
 where
 \be
 a_{\alpha ij}=\int_{\om}x_{\alpha}\lan\nabla u_i, \nabla u_j\ra=a_{\alpha ji}.
\en
We have
\be
\phi_{\alpha i}|_{\pa \om}=\left.\fr{\pa \phi_{\alpha i}}{\pa \nu}\right|_{\pa \om}=\cdots\left. \fr{\pa^{l-1} \phi_{\alpha i}}{\pa \nu^{l-1}}\right|_{\pa \om}=   0,
\en
\be
(\phi_{\alpha i}, u_j)_D= \int_{\om}\langle\nabla \phi_{\alpha i}, \nabla u_j\rangle=0, \ \ \forall j=1,\cdots, k.
\en
It then follows from the Rayleigh-Ritz inequality for $\laa_{k+1}$ that
\be
\laa_{k+1}\int_{\om}|\nabla \phi_{\alpha i}|^2\leq \int_\Omega\phi_{\alpha i}(-\Delta)^l\phi_{\alpha i}, \ \ \forall\alpha =1,\cdots, n, \ \ i=1,\cdots, k.
\en
After some  calculations, we have (cf. (2.36) in \cite{[JLWX1]})
\be
\int_{\om}\phi_{\alpha i}(-\Delta)^l\phi_{\alpha i}
&=&\int_{\om}(-1)^l\left\{(-l+1)u_i\D^{l-1}u_i+(2l^2-4l+3)(\D^{l-2}u_i)_{,\alpha}u_{i,\alpha}\right\}
\\ \no & & +\laa_i\left\{\int_{\om}x_{\alpha}^2|\na u_i|^2- \int_{\om}u_i^2\right\}-\sum_{j=1}^k\laa_ja_{\alpha ij}^2.
\en
It is easy to see that
\be
||x_{\alpha}\nabla u_i||^2=||\nabla h_{\alpha i}||^2+||W_{\alpha i}||^2, \ \ ||\nabla h_{\alpha i}||^2=||\nabla \phi_{\alpha i}||^2+\sum_{j=1}^ka_{\alpha ij}^2,
\en
where for a vector field $Z$ on $\om$, $||Z||^2=\int_{\om} |Z|^2$.
Combining (2.14)-(2.16), we infer
\be
(\laa_{k+1}-\laa_i)
||\nabla \phi_{\alpha i}||^2&\leq& \int_{\om}(-1)^l\left\{(-l+1)u_i\D^{l-1}u_i+(2l^2-4l+3)(\D^{l-2}u_i)_{,\alpha}u_{i,\alpha}\right\}
\\ \no & & -\laa_i(||u_i||^2-||W_{\alpha i}||^2)+\sum_{j=1}^k(\laa_i-\laa_j)a_{\alpha ij}^2,
\en
Observe that $\na(x_{\alpha}u_i)=u_i\na x_{\alpha}+ x_{\alpha}\na u_i\in {\bf H}_{1, D}^2(\om )$. Set  $y_{\alpha i}=x_{\alpha} u_i-h_{\alpha i}$; then
\be\no
u_i\na x_{\alpha}=\na y_{\alpha i}-W_{\alpha i}.
\en
and so
\be
||u_{i}||^2=||u_i \na x_{\alpha}||^2=||W_{\alpha i}||^2+||\na y_{\alpha i}||^2.
\en
Substituting (2.18) into (2.17), we get
\be\no & &
(\laa_{k+1}-\laa_i)
||\nabla \phi_{\alpha i}||^2\\ \no &\leq& \int_{\om}(-1)^l\left\{(-l+1)u_i\D^{l-1}u_i+(2l^2-4l+3)(\D^{l-2}u_i)_{,\alpha}u_{i,\alpha}\right\}
\\ \no & & -\laa_i ||\na y_{\al i}||^2+\sum_{j=1}^k(\laa_i-\laa_j)a_{\alpha ij}^2.
\en
Summing on $\alpha$ from $1$ to $n$, we have
\be& &
(\laa_{k+1}-\laa_i)\sum_{\alpha =1}^n
||\nabla \phi_{\alpha i}||^2\\ \no &\leq& \int_{\om}(-1)^l\left\{n(-l+1)u_i\D^{l-1}u_i+(2l^2-4l+3)\lan \na(\D^{l-2}u_i),\na u_i\ra\right\}
\\ \no & & -\laa_i\sum_{\alpha=1}^n ||\na y_{\al i}||^2+\sum_{\alpha=1}^n\sum_{j=1}^k(\laa_i-\laa_j)a_{\alpha ij}^2
\\ \no &=& (2l^2+(n-4)l+3-n)\int_{\om}u_i(-\D)^{l-1}u_i-\laa_i\sum_{\alpha=1}^n ||\na y_{\al i}||^2+\sum_{\alpha=1}^n\sum_{j=1}^k(\laa_i-\laa_j)a_{\alpha ij}^2.
\en
Using   the divergence  theorem,    one can show    that    (cf. \cite{[CY3]}, \cite{[JLWX1]})
\be
-2\int_{\om} x_{\alpha}\lan \na u_i, \na \lan \nabla u_i, \na x_{\alpha}\ra\ra =1.
\en
Set
\be\no
d_{\alpha ij}=\int_{\om}\lan \na \lan \nabla u_i, \na x_{\alpha}\ra, \na u_j\ra ,
\en
then $d_{\alpha ij}=-d_{\alpha ji}$ and we have from (2.7), (2.8),  (2.10) and (2.20) that
\be\no
1&=&-2\int_{\om} \lan \na h_{\alpha i}, \na \lan \nabla u_i, \na x_{\alpha}\ra\ra
\\ \no &=& -2\int_{\om} \lan \na \phi_{\alpha i}, \na \lan \nabla u_i, \na x_{\alpha}\ra\ra -2\sum_{j=1}^k a_{\alpha ij}d_{\alpha ij}.
\en
Thus, we have
\be
& & (\laa_{k+1}-\laa_i)^2\left(1+2\sum_{j=1}^k a_{\alpha ij}d_{\alpha ij}\right)\\ \no
&=&(\laa_{k+1}-\laa_i)^2\left( -2\na \phi_{\alpha i}, \na u_{i,\alpha}-\sum_{j=1}^k d_{\alpha ij}\na u_j\right)\\ \no &\leq&\delta_i (\laa_{k+1}-\laa_i)^3||\na \phi_{\alpha i}||^2+\fr 1{\delta_i} (\laa_{k+1}-\laa_i)\left(||\na u_{i,\alpha}||^2-\sum_{j=1}^k d_{\alpha ij}^2\right),
\en
where $u_{i,\alpha}= \lan\na u_i, \na x_{\alpha}\ra$.
Summing on $\alpha$ from $1$ to $n$, we have by using (2.19) that
\be\no
& & (\laa_{k+1}-\laa_i)^2\left(n+2\sum_{\alpha=1}^n\sum_{j=1}^k a_{\alpha ij}d_{\alpha ij}\right)\\ \no &\leq&
\delta_i (\laa_{k+1}-\laa_i)^2\left((2l^2+(n-4)l+3-n)\int_{\om}u_i(-\D)^{l-1}u_i-\laa_i\sum_{\alpha=1}^n ||\na y_{\al i}||^2\right.\\ \no & &\left.+\sum_{\alpha=1}^n\sum_{j=1}^k(\laa_i-\laa_j)a_{\alpha ij}^2\right)+\fr 1{\delta_i} (\laa_{k+1}-\laa_i)\left(\fa ||\na u_{i,\alpha}||^2-\fa\sum_{j=1}^k d_{\alpha ij}^2\right).
\en
Summing on $i$ from $1$ to $k$  and noticing the fact that $a_{\alpha ij}=a_{\alpha ji},\ d_{\alpha ij}=-d_{\alpha ji}$, one    gets
\be
& & n\sum_{i=1}^k(\laa_{k+1}-\laa_i)^2-2\fa\sum_{i, j=1}^k (\laa_{k+1}-\laa_i)(\laa_i-\laa_j)a_{\alpha ij}d_{\alpha ij}\\ \no
&\leq&\sum_{i=1}^k\delta_i (\laa_{k+1}-\laa_i)^2
\left((2l^2+(n-4)l+3-n)\int_{\om}u_i(-\D)^{l-1}u_i-\laa_i\sum_{\alpha=1}^n ||\na y_{\al i}||^2\right)\\ \no & &-\sum_{\alpha=1}^n\sum_{i, j=1}^k\delta_i(\laa_{k+1}-\laa_i)(\laa_i-\laa_j)^2a_{\alpha ij}^2
-\fa\sum_{i,j=1}^k \fr 1{\delta_i}(\laa_{k+1}-\laa_i)d_{\alpha ij}^2\\ \no & & +\sum_{i=1}^n\fr 1{\delta_i} (\laa_{k+1}-\laa_i)\fa ||\na u_{i,\alpha}||^2
\\ \no & &+\sum_{\alpha=1}^n\sum_{i, j=1}^k\delta_i(\laa_{k+1}-\laa_i)(\laa_i-\laa_j)^2a_{\alpha ij}^2
+\fa\sum_{i,j=1}^k \delta_i(\laa_{k+1}-\laa_i)^2(\laa_i-\laa_j)a_{\alpha ij}^2.
\en
Since  $\{\delta_i\}_{i=1}^k$ is a non-increasing monotone sequence, we have
\be\no & &
\sum_{\alpha=1}^n\sum_{i, j=1}^k\delta_i(\laa_{k+1}-\laa_i)(\laa_i-\laa_j)^2a_{\alpha ij}^2
+\fa\sum_{i,j=1}^k \delta_i(\laa_{k+1}-\laa_i)^2(\laa_i-\laa_j)a_{\alpha ij}^2
\\ \no &=&
\fr 12 \fa \sum_{i, j=1}^k(\laa_{k+1}-\laa_i)(\laa_{k+1}-\laa_j)(\laa_i-\laa_j)(\delta_i-\delta_j)a_{\alpha ij}^2\leq 0.
\en
We conclude from (2.22) that
\be
& & n\sum_{i=1}^k(\laa_{k+1}-\laa_i)^2\\ \no
&\leq&\sum_{i=1}^k\delta_i (\laa_{k+1}-\laa_i)^2
\left((2l^2+(n-4)l+3-n)\int_{\om}u_i(-\D)^{l-1}u_i-\laa_i\sum_{\alpha=1}^n ||\na y_{\al i}||^2\right)\\ \no & &\\ \no & & +\sum_{i=1}^n\fr 1{\delta_i} (\laa_{k+1}-\laa_i)\fa ||\na u_{i,\alpha}||^2.
\en
It  follows from    the divergence  theorem and Lemma 2.1   that
\be\no
\sum_{\alpha =1}^k||\na u_{i,\alpha}||^2&=&-\int_{\om}\sum_{\alpha =1}^k u_{i,\alpha}\D  u_{i,\alpha}\\ \no
&=&-\int_{\om}\sum_{\alpha =1}^k u_{i,\alpha}(\D  u_i)_{,\alpha}
\\ \no &=&
\int_{\om}\sum_{\alpha =1}^k u_{i,\alpha\alpha }\D  u_i
\\ \no &=&
\int_{\om}(\D  u_{i})^2 \\ \no &=&
\int_{\om}u_i\D^2  u_{i}\\ \no &\leq&\laa_i^{1/(l-1)},
\en
where $u_{i,\alpha\alpha }=\fr{\pa^2 u_i}{\pa x_{\alpha}^2}$.
Thus, we have
\be
& & n\sum_{i=1}^k(\laa_{k+1}-\laa_i)^2
\\ \no
&\leq   &\sum_{i=1}^k \delta_i (\laa_{k+1}-\laa_i)^2\left( (2l^2+(n-4)l+3-n) \int_{\om}u_i(-\D)^{l-1}u_i-\sum_{\al=1}^n\laa_i    ||\na y_{\alpha i}||^2
\right)
\\ \no & &
+ \sum_{i=1}^k \fr 1{\delta_i} (\laa_{k+1}-\laa_i)\laa_i^{1/(l-1)}.
\en
\vskip0.3cm
\noindent
Before we can finish the proof of Theorem  1.2, we shall need  two lemmas.
\vskip0.3cm
\noindent
{\bf Lemma 2.2.} For any $i$, we have
\be (n-2l-2)\int_{\om}  u_i(-\Delta)^{l-1}u_i = n\laa_i ||u_i||^2-4\laa_i\sum_{\alpha=1}^n ||\na y_{\alpha i}||^2.
\en
{\it Proof.} When   $l=2$,  the above   formula has been    proved  by  Cheng-Yang  in  \cite{[CY8]}. We  only    consider    the case    that
 $l>2$. In this case,  we  conclude    from    the boundary    condition   on  $u_i$   that   $ y_{\al i}|_{\partial \om}=\na  y_{\al i}|_{\partial \om}=\D y_{\al i}|_{\partial \om}=0$. Using the divergence theorem, we have
\be \int_{\om} x_{\al} u_i\lan\na x_{\al}, \na(\D^{l-1}u_i)\ra&=&\int_{\om} x_{\al} u_i\D^{l-1}\lan\na x_{\al},\na u_i\ra\\ \no &=&
\int_{\om} \D^{l-1}(x_{\al} u_i)\lan\na x_{\al},\na u_i\ra\\ \no &=&
-\int_{\om} \lan u_i\na x_{\al},\na(\D^{l-1}(x_{\al} u_i))\ra\\ \no &=&
-\int_{\om} \lan \na y_{\al i},\na(\D^{l-1}(x_{\al} u_i))\ra\\ \no &=&
\int_{\om}  y_{\al i}\D^{l}(x_{\al} u_i)\\ \no &=&
\int_{\om}  y_{\al i}(2l\lan\na(\D^{l-1}u_i), \na x_{\al}\ra+x_{\al}\D^lu_i)\\ \no &=&
-2l\int_{\om}\D^{l-1}u_i\lan \na y_{\al i}, \na x_{\al}\ra+\laa_i (-1)^{l-1}\int_{\om}y_{\al i}x_{\al}\D u_i,
\en
\be\int_{\om}y_{\al i}x_{\al}\D u_i&=&-\int_{\om}\lan\na y_{\al i},\ x_{\al}\na u_i\ra
-\int_{\om}y_{\al i}\lan\na x_{\al}, \na u_i\ra\\ \no &=&-\int_{\om}\lan\na y_{\al i},\ x_{\al}\na u_i\ra
+\int_{\om}\lan\na y_{\al i},\ u_i\na x_{\al}\ra\\ \no &=&-\int_{\om}\lan\na y_{\al i},\ x_{\al}\na u_i\ra
+||\na y_{\al i}||^2,
\en
\be
\int_{\om}\lan\na y_{\al i},\ x_{\al}\na u_i\ra&=&
\int_{\om}\lan\na y_{\al i},\na h_{\al i}\ra\\ \no &=&
\int_{\om}\lan\na y_{\al i},\na(x_{\al}u_i)-\na y_{\al i}\ra\\ \no &=&
\int_{\om}\lan\na y_{\al i},\na(x_{\al}u_i)\ra-||\na y_{\al i}||^2 \\ \no &=&
\int_{\om}\lan u_i\na x_{\al},\na(x_{\al}u_i)\ra-||\na y_{\al i}||^2 \\ \no &=&
||u_i||^2+\int_{\om}\lan u_i\na x_{\al},\ x_{\al}\na u_i\ra-||\na y_{\al i}||^2
\\ \no &=&
||u_i||^2-\fr 14\int_{\om}u_i^2\D x_{\al}^2-||\na y_{\al i}||^2 \\ \no &=&\fr 12 ||u_i||^2-||\na y_{\al i}||^2
\en
and
\be
\int_{\om}\D^{l-1}u_i\lan \na y_{\al i}, \na x_{\al}\ra &=&-\int_{\om}y_{\al i}\lan \na (\D^{l-1}u_i), \na x_{\al}\ra
\\ \no &=&-\int_{\om}y_{\al i}  \D\lan\na (\D^{l-2} u_i), \na x_{\al}\ra
\\ \no&=&
\int_{\om}\lan\na  y_{\al i}, \na\lan \na(\D^{l-2}u_i), \na x_{\al}\ra\ra
\\ \no&=&
\int_{\om}\lan u_i\na  x_{\al }, \na\lan \na(\D^{l-2}u_i), \na x_{\al}\ra\ra
\\ \no&=&
-\int_{\om}\lan \na  x_{\al },  \na(\D^{l-2}u_i)\ra \lan\na u_i, \na x_{\al}\ra.
\en
It follows from (2.26)-(2.29) that
\be & &\int_{\om} x_{\al} u_i\lan\na x_{\al}, \na(\D^{l-1}u_i)\ra\\ \no &=&
2l\int_{\om}\lan \na  x_{\al },  \na(\D^{l-2}u_i)\ra \lan\na u_i, \na x_{\al}\ra
+(-1)^{l-1}\laa_i\left(-\fr 12||u_i||^2+2||\na y_{\al i}||^2\right).
\en
Since
\be\no
\D^{l-1}(x_{\alpha}u_i)=2(l-1)\lan \na (\D^{l-2}u_i), \na x_{\alpha}\ra +x_{\alpha}\D^{l-1} u_i,
\en
we get
\be
\int_{\om}x_{\alpha}u_i \lan\na x_{\alpha}, \na(\D^{l-1}u_i) \ra  &=&\int_{\om}x_{\alpha}u_i \D^{l-1}\lan \na u_i, \na x_{\alpha}\ra\\ \no &=&\int_{\om}\D^{l-1}(x_{\alpha}u_i) \lan \na u_i, \na x_{\alpha}\ra\\ \no &=&\int_{\om}\left(2(l-1)\lan \na (\D^{l-2}u_i), \na x_{\alpha}\ra +x_{\alpha}\D^{l-1} u_i\right) \lan \na u_i, \na x_{\alpha}\ra.
\en
On the other hand, we have
\be
\int_{\om}x_{\alpha}u_i \lan \na x_{\alpha}, \na(\D^{l-1}u_i)\ra = -\int_{\om}\D^{l-1}u_i(u_i+ x_{\alpha}\lan \na u_i, \na x_{\alpha}\ra).
\en
We obtain from (2.31) and (2.32) that
\be & &
\int_{\om}x_{\alpha}u_i \lan\na x_{\alpha}, \na(\D^{l-1}u_i)\ra
 \\ \no &=&
 \int_M\left\{(l-1)(\lan \na(\D^{l-2}u_i), \na x_{\alpha}\ra \lan \na u_i, \na x_{\alpha}\ra-\fr 12u_i\D^{l-1}u_i \right\}.
\en
Combining (2.30) and (2.33), we infer
 \be& &
 \int_M\left\{(l-1)\lan\na x_{\al}, \na(\D^{l-2}u_i)\ra \lan\na x_{\al},    \na u_i\ra-\fr 12u_i\D^{l-1}u_i \right\}
\\ \no &=&
2l\int_{\om}\lan \na  x_{\al },  \na(\D^{l-2}u_i)\ra \lan\na u_i, \na x_{\al}\ra
+(-1)^{l-1}\laa_i\left(-\fr 12||u_i||^2+2||\na y_{\al i}||^2\right).
\en
Summing on $\al$, we get (2.25).
\vskip0.3cm
\noindent
{\bf Lemma 2.3.}  For any $i$, we have
\be
\sum_{\al =1}^n ||W_{\al i}||^2 \geq \fr{n-1}{\laa_i^{1/(l-1)}}.
\en
{\it Proof.} Using the definition of $W_{\al i}$ and the    divergence  theorem and noticing (2.20), we have
\be& & \int_{\om}  \lan \na x_{\al}, W_{\al i}\ra\D u_i\\ \no &=& -\int_{\om} \lan \na u_i, \na\lan \na x_{\al}, W_{\al i}\ra\ra\\ \no
&=& -\int_{\om} \lan \na u_i, \na\left(\lan  x_{\al}\na u_i-\na h_{\al i}, \na x_{\al}\ra\right)\ra\\ \no
&=& -||\lan\na u_i, \na x_{\al}\ra||^2-\int_{\om} x_{\al}\lan \na u_i, \na \lan \na u_i, \na x_{\al}\ra\ra+\int_{\om}\lan \na u_i,
\na \lan \na h_{\al i}, \na x_{\al}\ra\ra\\ \no
&=& -||\lan\na u_i, \na x_{\al}\ra||^2-\int_{\om} x_{\al}\lan \na u_i,  \na \lan \na u_i, \na x_{\al}\ra\ra -\int_{\om} u_i\D
\lan \na h_{\al i}, \na x_{\al}\ra \\ \no
\\ \no
&=& -||\lan\na u_i, \na x_{\al}\ra||^2-\int_{\om} x_{\al}\lan \na u_i,  \na \lan \na u_i, \na x_{\al}\ra\ra-\int_{\om} u_i\lan \na(\D  h_{\al i}), \na x_{\al}\ra \\ \no
\\ \no
&=& -||\lan\na u_i, \na x_{\al}\ra||^2-\int_{\om} \lan x_{\al}\na u_i,  \na \lan \na u_i, \na x_{\al}\ra\ra -\int_{\om} u_i
\lan\na (\lan \na x_{\al}, \na u_i\ra +x_{\al} \D u_i), \na x_{\al}\ra
\\ \no
&=& -||\lan\na u_i, \na x_{\al}\ra||^2-\int_{\om} \lan \na (x_{\al} u_i),  \na \lan \na u_i, \na x_{\al}\ra\ra -\int_{\om} u_i
\lan\na (x_{\al} \D u_i), \na x_{\al}\ra
\\ \no
&=& -||\lan\na u_i, \na x_{\al}\ra||^2+\int_{\om}   \lan \na u_i, \na x_{\al}\ra\D (x_{\al} u_i) +\int_{\om}
x_{\al} \D u_i\lan\na  u_i, \na x_{\al}\ra
\\ \no
&=& ||\lan\na u_i, \na x_{\al}\ra||^2+2\int_{\om}
x_{\al} \D u_i\lan\na  u_i, \na x_{\al}\ra
\\ \no
&=& ||\lan\na u_i, \na x_{\al}\ra||^2-2\int_{\om}
 \lan\na u_i, \na (x_{\al}\lan\na  u_i, \na x_{\al}\ra)\ra
\\ \no
&=& -||\lan\na u_i, \na x_{\al}\ra||^2-2\int_{\om}
 x_{\al}\lan\na u_i, \na \lan\na  u_i, \na x_{\al}\ra\ra
\\ \no
&=& -||\lan\na u_i, \na x_{\al}\ra||^2+1.
\en
On the other hand,  for $\epsilon >0$, we have
\be
\int_{\om}\lan \na x_{\al}, W_{\al i}\ra\D u_i &=& \int_{\om}\lan \D u_i\na x_{\al}-\na \lan\na  u_i, \na x_{\al}\ra, W_{\al i}\ra
\\ \no &\leq& \fr{\epsilon}2||W_{\al i}||^2+\fr 1{2\epsilon}||\D u_i\na x_{\al}-\na \lan\na  u_i, \na x_{\al}\ra||^2.
\en
From (2.36), we have
\be
\sum_{\al=1}^n \int_{\om}\lan \na x_{\al}, W_{\al i}\ra\D u_i=n-1.
 \en
 Also, one can  check that
 \be
 \sum_{\al=1}^n||\D u_i\na x_{\al}-\na \lan\na  u_i, \na x_{\al}\ra||^2= (n-1)\int_{\om} u_i\D^2 u_i\leq(n-1)\laa_i^{1/(l-1)} .
 \en
 Thus we have from (2.37)-(2.39) that
 \be
 n-1\leq \fr{\epsilon}2\sum_{\al=1}^n||W_{\al i}||^2+\fr{n-1}{2\epsilon}\laa_i^{1/(l-1)}.
 \en
 Taking
 $$\epsilon =\sqrt{\fr{(n-1)\laa_i^{1/(l-1)}}{\sum_{\al=1}^n||W_{\al i}||^2}},$$
 we get (2.35). This completes the proof of Lemma 2.3.
 
\vskip0.3cm
\noindent
Let us  continue  the proof   of  Theorem 1.1.
 Since $||u_i||^2=||W_{\al i}||^2+||\na y_{\al i}||^2$, we have from (2.35) that
 \be
 n\laa_i ||u_i||^2&=&\laa_i\sum_{\al=1}^n||W_{\al i}||^2+\laa_i\sum_{\al=1}^n||\na y_{\al i}||^2\\ \no
 &\geq&(n-1)\laa_i^{(l-2)/(l-1)}+\laa_i\sum_{\al=1}^n||\na y_{\al i}||^2,
 \en
 which, combining with (2.25), implies that
 \be
 -\laa_i\sum_{\al=1}^n||\na y_{\al i}||^2\leq\fr{(n-2l-2)}3\int_{\om}u_i(-\D)^{l-1}u_i-\fr{(n-1)}3\laa_i^{(l-2)/(l-1)}.
 \en
 Substituting (2.42) into (2.24) and using Lemma 2.1, we get
 \be\no
& & n\sum_{i=1}^k(\laa_{k+1}-\laa_i)^2
\\ \no
&\leq   &\sum_{i=1}^k \delta_i(\laa_{k+1}-\laa_i)^2\left( \left(2l^2+(n-4)l+3-n+\fr{n-2l-2}3\right) \int_{\om}u_i(-\D)^{l-1}u_i
-\fr{(n-1)}3\laa_i^{(l-2)/(l-1)}
\right)
\\ \no & &
+ \sum_{i=1}^k \fr 1{\delta_i} (\laa_{k+1}-\laa_i)\sum_{\alpha =1}^n\laa_i^{1/(l-1)}\\ \no &\leq &
\sum_{i=1}^k\delta_i (\laa_{k+1}-\laa_i)^2\left( \left(2l^2+(n-4)l+3-n+\fr{n-2l-2}3-\fr{(n-1)}3\right) \laa_i^{(l-2)/(l-1)}
\right)
\\ \no & &
+\sum_{i=1}^k \fr 1{\delta_i}  (\laa_{k+1}-\laa_i)\laa_i^{1/(l-1)}
\\ \no & =&\sum_{i=1}^k \delta_i(\laa_{k+1}-\laa_i)^2\left(2l^2+\left(n-\fr{14}3\right)l+\fr 83-n\right)\laa_i^{(l-2)/(l-1)}\\ \no & &
+ \sum_{i=1}^k \fr 1{\delta_i} (\laa_{k+1}-\laa_i)\laa_i^{1/(l-1)}.
\en
 This completes the proof of Theorem 1.1.
\vskip0.3cm
\noindent
{\it Proof of Corollary 1.2.} By induction, one can show that
\be\no & &
\left\{ \sum_{i=1}^k  (\laa_{k+1}-\laa_i)^2\laa_i^{(l-2)/(l-1)} \right\} \left\{\sum_{i=1}^k (\laa_{k+1}-\laa_i)\laa_i^{1/(l-1)}\right\} \\ \no &\leq &
\left\{ \sum_{i=1}^k  (\laa_{k+1}-\laa_i)^2 \right\} \left\{\sum_{i=1}^k (\laa_{k+1}-\laa_i)\laa_i\right\},
\en
which, combining with (1.12), gives (1.13).

\vskip0.3cm
\noindent
{\it Proof of Theorem 1.2.} We use  the same notations as in the beginning of this section and take $M$ to be the unit $n$-sphere $S^n$.
Let $x_1, x_2,\cdots, x_{n+1}$ be the standard  coordinate functons of the Euclidean space  ${\bf R}^{n+1}$, then
$$S=\left\{(x_1,\dots, x_{n+1})\in  {\bf R}^{n+1}; \sum_{\alpha=1}^{n+1}x_{\alpha}^2=1\right\}.$$ It is well known that
\be
\Delta x_{\alpha}=-nx_{\alpha},\ \ \  \alpha =1,\cdots, n+1.
\en
As in the proof of Theorem 1.1, we
 decompose the vector-valued functions $x_{\alpha}\nabla u_i$ as
\be
x_{\alpha}\nabla u_i=\nabla h_{\alpha i}+ W_{\alpha i},
\en
where
$h_{\alpha i}\in H_{l,D}^2(\om),$ $\nabla h_{\alpha i}$ is the projection of $x_{\alpha} \nabla u_i$ in ${\bf H}_{1, D}^2(\om)$, $W_{\alpha i}\ \bot\ {\bf H}_{1, D}^2(\om)$
 and
\be
W_{\alpha i}|_{\pa\om}\ =\ 0, \ \  {\rm div\ } W_{\alpha i}=0.
\en
We also consider the  functions $\phi_{\alpha i}: \om\rightarrow {\bf R}$, given by
\be
 \phi_{\alpha i}=h_{\alpha i}-\sum_{j=1}^kb_{\alpha ij}u_j, \ \
 b_{\alpha ij}=\int_{\om}x_{\alpha}\lan\nabla u_i, \nabla u_j\ra=b_{\alpha ji}.
\en
 Then
\be\no
\phi_{\alpha i}|_{\pa \om}=\left.\fr{\pa \phi_{\alpha i}}{\pa \nu}\right|_{\pa \om}=\cdots\left. \fr{\pa^{l-1} \phi_{\alpha i}}{\pa \nu^{l-1}}\right|_{\pa \om}=   0,
\en
\be\no
(\phi_{\alpha i}, u_j)_D= \int_{\om}\langle\nabla \phi_{\alpha i}, \nabla u_j\rangle=0, \ \ \forall j=1,\cdots, k
\en
and we have the basic Rayleigh-Ritz inequality for $\laa_{k+1}$ :
\be
\laa_{k+1}\int_{\om}|\nabla \phi_{\alpha i}|^2\leq \int_D\phi_{\alpha i}(-\Delta)^l\phi_{\alpha i}, \ \ \forall\alpha =1,\cdots, n, \ \ i=1,\cdots, k.
\en
We have
\be
\Delta \phi_{\alpha i}= \langle \nabla x_{\alpha}, \nabla u_i\rangle+x_{\alpha}\Delta u_i-\sum_{j=1}^kb_{\alpha ij}\Delta u_j
\en
and from (2.56) in \cite{[JLWX1]},
\be
& & \int_{\om}\phi_{\alpha i}(-\Delta)^l\phi_{\alpha i}\\ \no & &
=\int_{\om}(-1)^l(\langle \nabla x_{\alpha}, \nabla u_i\rangle +x_{\alpha}\Delta u_i)\Delta^{l-2}(\langle \nabla x_{\alpha}, \nabla u_i\rangle+x_{\alpha}\Delta u_i)
-\sum_{j=1}^k\laa_jb_{\alpha ij}^2.
\en
For a function $g$ on $\om$, we have (cf. (2.31) in \cite{[WX2]})
\be
\Delta \langle \nabla x_{\alpha}, \nabla g\rangle = -2x_{\alpha}\Delta g+\langle \nabla x_{\alpha}, \nabla((\Delta +n-2)g)\rangle.
\en
For each $q=0, 1,\cdots$, thanks to (2.43) and (2.50), there are polynomials
$F_q$ and $G_q$ of degree $q$ such that
\be
\Delta^{q}(\langle \nabla x_{\alpha}, \nabla u_i\rangle+x_{\alpha}\Delta u_i)
=x_{\alpha}F_q(\D)\D u_i+\lan \na x_{\alpha}, \na (G_q(\D)u_i)\ra.
\en
It is obvious that
\be
F_0=1, \ \ G_0=1.
\en
It follows from (2.43) and (2.50) that
\be
\D (x_{\alpha}\Delta u_i+\langle \nabla x_{\alpha}, \nabla u_i\rangle)=x_{\alpha}(\D-(n+2))\D u_i+\langle \nabla x_{\alpha}, \na((3\D+n-2)u_i)\ra,
\en
which gives
\be
F_1(t)=t-(n+2), \ \ G_1(t)=3t+n-2.
\en
Also, when $q\geq 2$, we have (cf. (2.65) and (2.66) in \cite{[JLWX1]})
\be& &
 F_q(t)=(2t-2)F_{q-1}(t)-(t^2+2t-n(n-2))F_{q-2}(t), \ \ q=2,\cdots, \\  & & G_q(t)=(2t-2)G_{q-1}(t)-(t^2+2t-n(n-2))G_{q-2}(t), \ \ q=2,\cdots.
\en
\vskip0.3cm
For each $q=1, 2, \cdots, $ let us set
 $$\Phi_q(t)=tF_{q-1}(t)-G_{q-1}(t).$$ We conclude from (2.52), (2.54)-(2.56) that the polynomials
$\Phi_q, \ q=1, 2, \cdots, $ are defined inductively by (1.15) and (1.16).  Substituting
 \be
\Delta^{l-2}(\langle \nabla x_{\alpha}, \nabla u_i\rangle+x_{\alpha}\Delta u_i)
=x_{\alpha}F_{l-2}(\D)\D u_i+\lan \na x_{\alpha}, \na (G_{l-2}(\D)u_i)\ra
\en
into (2.49), we get
\be
& &\int_{\om}\phi_{\alpha i}(-\Delta)^l\phi_{\alpha i}\\ \no &=&
\int_{\om}(-1)^l(\langle \nabla x_{\alpha}, \nabla u_i\rangle \lan \na x_{\alpha}, \na (G_{l-2}(\D)u_i)\ra
+\lan x_{\alpha}\na x_{\alpha}, \Delta u_i\na (G_{l-2}(\D)u_i)+(F_{l-2}(\D)\D u_i) \na u_i\ra)\\ \no & & +\int_{\om}(-1)^l x_{\alpha}^2\Delta u_iF_{l-2}(\D)(\D u_i) -\sum_{j=1}^k\laa_jb_{\alpha ij}^2.
\en
Summing over $\alpha $  and noticing
\be
 \sum_{\alpha=1}^{n+1}x_{\alpha}^2=1, \ \  \sum_{\alpha=1}^{n+1}\langle \nabla x_{\alpha}, \nabla u_i\rangle \lan \na x_{\alpha}, \na (G_{l-2}(\D)u_i)\ra = \langle  \nabla u_i , \na (G_{l-2}(\D)u_i)\ra,
\en
we infer
\be
& &\sum_{\alpha=1}^{n+1}\int_{\om}\phi_{\alpha i}(-\Delta)^l\phi_{\alpha i}\\ \no &=&
\int_{\om}(-1)^l\langle \nabla u_i, \na (G_{l-2}(\D)u_i)\ra
 +\int_{\om}(-1)^l \Delta u_iF_{l-2}(\D)(\D u_i) -\sum_{\alpha=1}^{n+1}\sum_{j=1}^k\laa_jb_{\alpha ij}^2
\\ \no &=&
\int_{\om}(-1)^{l-1} u_i \D (G_{l-2}(\D)u_i)
 +\int_{\om}(-1)^l  u_i\D (F_{l-2}(\D)(\D u_i)) -\sum_{\alpha=1}^{n+1}\sum_{j=1}^k\laa_jb_{\alpha ij}^2
\\ \no &=&
\int_{\om}(-1)^{l} u_i \left( \D F_{l-2}(\D)-G_{l-2}(\D)\right)(\D u_i)
-\sum_{\alpha=1}^{n+1}\sum_{j=1}^k\laa_jb_{\alpha ij}^2
\\ \no &=&
\int_{\om}(-1)^{l} u_i \Phi_{l-1}(\D)(\D u_i)
-\sum_{\alpha=1}^{n+1}\sum_{j=1}^k\laa_jb_{\alpha ij}^2
\\ \no &=&
\int_{\om}(-1)^{l} u_i \left( \D^{l-1}-a_{l-2}\D^{l-2}+\cdots +(-1)^{l-2}a_1 \D -(n-2)^{l-2}\right)(\D u_i)
-\sum_{\alpha=1}^{n+1}\sum_{j=1}^k\laa_jb_{\alpha ij}^2
\\ \no &= &
\laa_i+ (-1)^l(n-2)^{l-2}+\sum_{j=1}^{l-2}a_j\int_{\om} u_i (-\D)^{j+1}u_i-\sum_{\alpha=1}^{n+1}\sum_{j=1}^k\laa_jb_{\alpha ij}^2.
\en
Set
\be
H_i=   (-1)^l(n-2)^{l-2}+\sum_{j=1}^{l-2}a_j^{+}\laa_i^{j/(l-1)},
\en
then it is  easy to check from   Lemma 2.1  that
\be
(-1)^l(n-2)^{l-2}+\sum_{j=0}^{l-2}a_j\int_{\om} u_i (-\D)^{j+1}u_i\leq H_i.
\en
Substituting (2.62) into (2.60), we have
\be\sum_{\alpha=1}^{n+1}\int_{\om}\phi_{\alpha i}(-\Delta)^l\phi_{\alpha i}
\leq   \laa_i+ H_i-\sum_{\alpha=1}^{n+1}\sum_{j=1}^k\laa_jb_{\alpha ij}^2.
\en
Observe from (2.44) and (2.46) that
\be
||x_{\alpha}\na u_i||^2=||\na h_{\alpha i}||^2+||W_{\alpha i}||^2=
||\na \phi_{\alpha i}||^2+||W_{\alpha i}||^2+\sum_{j=1}^kb_{\alpha ij}^2.
\en
Summing over $\alpha$, one gets
\be
1=\sum_{\alpha=1}^{n+1}\left(||\na \phi_{\alpha i}||^2+||W_{\alpha i}||^2+\sum_{j=1}^kb_{\alpha ij}^2\right).
\en
Combining  (2.47), (2.63) and (2.65), we get
\be
\sum_{\alpha=1}^{n+1}(\laa_{k+1}-\laa_i)||\nabla \phi_{\alpha i}||^2\leq H_i
+\sum_{\alpha=1}^{n+1}\laa_i||W_{\alpha i}||^2+\sum_{\alpha=1}^{n+1}\sum_{j=1}^k(\laa_i-\laa_j)b_{\alpha ij}^2.
\en
Set
\be
Z_{\alpha i}=\nabla\langle \nabla x_{\alpha},\ \nabla u_i\rangle-\fr{n-2}2x_{\alpha}\nabla u_i,\ \
c_{\alpha ij}=\int_{\om}\langle\nabla u_j,\ Z_{\alpha i}\rangle
;
\en
then
$c_{\alpha ij}=-c_{\alpha ji}$ (cf. Lemma in \cite{[WX2]}). By using the same arguments as in the proof of (2.37) in \cite{[WX2]}, we have
\be & &
(\laa_{k+1}-\laa_i)^2\left(2||\langle\nabla x_{\alpha}, \nabla u_i\rangle||^2 +\int_{\om} \left\langle\nabla x_{\alpha}^2,\ \Delta u_i\nabla u_i\right\rangle +(n-2)|| x_{\alpha}\nabla u_i||^2+2\sum_{j=1}^k b_{\alpha ij}c_{\alpha ij}\right)
\\ \no  &\leq &\delta_i (\laa_{k+1}-\laa_i)^3||\nabla \phi_{\alpha i}||^2
+\fr{\laa_{k+1}-\laa_i}{\delta_i}\left(||Z_{\alpha i}||^2-\sum_{j=1}^kc_{\alpha ij}^2\right)
+(n-2)(\laa_{k+1}-\laa_i)^2 ||W_{\alpha i}||^2.
\en

Since
\be
\sum_{\alpha=1}^{n+1}||\langle\nabla x_{\alpha}, \nabla u_i\rangle||^2=\int_{\om} |\na u_i|^2=1,
\en
we have by
summing over $\alpha$ in (2.68) from 1 to $n+1$  that
\be & &
(\laa_{k+1}-\laa_i)^2\left(n+2\sum_{\alpha=1}^{n+1}\sum_{j=1}^k b_{\alpha ij}c_{\alpha ij}\right)
\\ \no  &\leq &\delta_i \sum_{\alpha=1}^{n+1}(\laa_{k+1}-\laa_i)^3||\nabla \phi_{\alpha i}||^2
+\fr{\laa_{k+1}-\laa_i}{\delta_i}\sum_{\alpha=1}^{n+1}\left(||Z_{\alpha i}||^2-\sum_{j=1}^kc_{\alpha ij}^2\right)
\\ \no & & +(n-2)\sum_{\alpha=1}^{n+1}(\laa_{k+1}-\laa_i)^2 ||W_{\alpha i}||^2.
\en
From (2.77), (2.78) and (2.80) in \cite{[JLWX1]}, we have
\be
\laa_i^{1/(l-1)}-(n-2)>0,
\en

\be & &
\sum_{\alpha=1}^{n+1}||Z_{\alpha i}||^2
 \leq\laa_i^{1/(l-1)}+\fr{(n-2)^2}4
\en
and
\be
\sum_{\alpha=1}^{n+1}||W_{\alpha i}||^2\leq 1 -\fr
1{\laa_i^{1/(l-1)}-(n-2)}. \en

We have by combining (2.66), (2.70), (2.72) and (2.73)   that
\be & &
(\laa_{k+1}-\laa_i)^2\left(n+2\sum_{\alpha=1}^{n+1}\sum_{j=1}^k b_{\alpha ij}c_{\alpha ij}\right)
\\ \no &\leq &
\delta_i (\laa_{k+1}-\laa_i)^2\left( H_i
+\sum_{\alpha=1}^{n+1}\sum_{j=1}^k(\laa_i-\laa_j)b_{\alpha ij}^2\right)
\\ \no & & +\fr{\laa_{k+1}-\laa_i}{\delta_i}\left(||Z_{\alpha i}||^2-\sum_{\alpha=1}^{n+1}\sum_{j=1}^kc_{\alpha ij}^2\right)
+\sum_{\alpha=1}^{n+1}(\laa_{k+1}-\laa_i)^2(\delta_i \laa_i+n-2) ||W_{\alpha i}||^2\\ \no
&\leq &
\delta_i (\laa_{k+1}-\laa_i)^2\left( H_i
+\sum_{\alpha=1}^{n+1}\sum_{j=1}^k(\laa_i-\laa_j)b_{\alpha ij}^2\right)
\\ \no & & +\fr{\laa_{k+1}-\laa_i}{\delta_i}\left(\left(\laa_i^{1/(l-1)}+\fr{(n-2)^2}4\right)-\sum_{\alpha=1}^{n+1}\sum_{j=1}^kc_{\alpha ij}^2\right)
\\ \no & & +(\laa_{k+1}-\laa_i)^2(\delta_i \laa_i+n-2) \left( 1
-\fr 1{\laa_i^{1/(l-1)}-(n-2)}\right).
\en
Since   $\{\delta_i\}_{i=1}^k$ is a positive non-increasing monotone sequence, we have
\be& &
2\sum_{\alpha=1}^{n+1}\sum_{i, j=1}^k(\laa_{k+1}-\laa_i)^2 b_{\alpha ij}c_{\alpha ij}\\ \no &\geq&
\sum_{\alpha=1}^{n+1}\sum_{i, j=1}^k \delta_i (\laa_{k+1}-\laa_i)^2 (\laa_i-\laa_j)b_{\alpha ij}^2-\sum_{\alpha=1}^{n+1}\sum_{i, j=1}^k\fr{\laa_{k+1}-\laa_i}{\delta_i} c_{\alpha ij}^2.
\en
Hence, by summing over $i$ from 1 to $k$ in (2.74), we infer
\be \no & &
n\sum_{i=1}^k(\laa_{k+1}-\laa_i)^2
\\ \no &\leq &
 \sum_{i=1}^k(\laa_{k+1}-\laa_i)^2 \left(\delta_i H_i+(\delta_i \laa_i+n-2)\left( 1
-\fr 1{\laa_i^{1/(l-1)}-(n-2)}\right)\right)\\ \no & &
 +\sum_{i=1}^k\fr{\laa_{k+1}-\laa_i}{\delta_i}\left(\laa_i^{1/(l-1)}+\fr{(n-2)^2}4\right).
\en
That is
\be  \no & &
\sum_{i=1}^k(\laa_{k+1}-\laa_i)^2\left(2+\fr{n-2}{\laa_i^{1/(l-1)}-(n-2)}\right)
\\ \no &\leq &
 \sum_{i=1}^k(\laa_{k+1}-\laa_i)^2\delta_i\left(H_i+\laa_i\left( 1
-\fr 1{\laa_i^{1/(l-1)}-(n-2)}\right)\right)\\ \no & &+
\sum_{i=1}^k\fr {(\laa_{k+1}-\laa_i)}{\delta_i}\left(\laa_i^{1/(l-1)}+\fr{(n-2)^2}4\right)\\ \no
&=&
\sum_{i=1}^k(\laa_{k+1}-\laa_i)^2\delta_i S_i+
\sum_{i=1}^k\fr {(\laa_{k+1}-\laa_i)}{\delta_i}\left(\laa_i^{1/(l-1)}+\fr{(n-2)^2}4\right),
\en
where $S_i$ is given by (1.19).
 This completes the proof of Theorem 1.2.

\vs
\noindent Qing-Ming Cheng (cheng@ms.saga-u.ac.jp) and Xuerong Qi (qixuerong609@gmail.com)

\noindent Department of Mathematics

\noindent Faculty of Science and Engineering

\noindent Saga University

\noindent Saga 840-8502, Japan

\vs
\noindent Qiaoling Wang (wang@mat.unb.br) and Changyu Xia (xia@mat.unb.br)

\noindent Departamento de Matem\'{a}tica

\noindent Universidade de Bras\'{\i}lia

\noindent 70910-900-Bras\'{\i}lia-DF, Brazil

\end{document}